\documentclass[a4paper,10pt]{article}

\def \Z {{\mathbf {Z}}}
\def \N {{\mathbf {N}}}
\def \F {{\cal F}}
\title{ Spectral multiplicity for  powers of  weakly mixing automorphisms }
\author{V.V. Ryzhikov\footnote{This work is partially  supported by  the grant
NSh 8508.2010.1.}}

\begin{document}
\Large
\date{11.02.2011}
\maketitle

\begin{abstract}{\large We study the behavior of  maximal   multiplicities
$mm (R^n)$ for the powers of a weakly mixing automorphism $R$.
For  some special  infinite set $A$ we show the existence of a weakly mixing 
rank-one automorphism $R$ such that  
$mm (R^n)=n$ and $mm(R^{n+1}) =1$
for all  $n\in A$. 
Moreover, the cardinality $cardm(R^n)$
of the set of spectral multiplicities for 
$R^n$ is not bounded.  We have $cardm(R^{n+1})=1$ and 
$cardm(R^n)=2^{m(n)}$,  $m(n)\to\infty$, $n\in A$.  
We also construct another weakly mixing automorphism $R$ with the following 
properties: $mm(R^{n}) =n$ for $n=1,2,3, \dots, 2009, 2010$ but $mm(T^{2011}) =1$,
all powers $(R^{n})$ have homogeneous spectrum, and the set of limit
points of the sequence $\{\frac{mm (R^n)}{n} : n\in \N \}$ is infinite. }
\end{abstract}

\section{ Introduction} 
We recall  the notion of  maximal 
spectral multiplicity  $mm (U)$ of a unitary operator~$U$.
Consider all representations of the  form
$$U\equiv W\oplus V\oplus V\oplus\dots \oplus V\dots\ .$$
The maximal number of such copies of $V$ is denoted by~$mm (U)$. 
   An automorphism $T$ acting  on a measure space $(X,\mu)$, $\mu(X)=1$,  induces 
a unitary operator $Uf=fT$ (acting on $L_2(X,\mu)$, $f\in L_2$). We write  $mm (T)$
instead of $mm (U)$.
The property of an automorphism to be weakly mixing means that its spectral measure
is continuous. 
Weak mixing can also be defined as the existence of a (mixing) sequence $m_i$
such that for all $f,g\in L_2(X,\mu)$
$$\int f U^{m_i}g  d\mu \to  \int f d\mu\int g d\mu.$$

In ergodic theory there are examples 
of weakly mixing transformations with the following properties:

(1) $mm(T^{n}) =1$,  $n=1,2,3,4,\dots$;

(2) $mm(T^{n}) =m$,  $n=1,2,3,4,\dots$;

(3) $mm(T^{n}) =2mn$, $n=1,2,3,4,\dots$.

For example, the latter follows from the results of Ageev [Ag1] and Lemanczyk [L]
on realization of a Lebesgue component of multiplicity~2 (here $mm(T) =2$ implies that
$mm(T^{2}) =4$, $mm(T^{3}) =6, \dots$).
Property (1) is known to hold for  

Gaussian systems with simple spectrum,  

Ornstein's stochastic constructions [Ab], 

some mixing staircase constructions [R].  
\\
The aim of this paper is to give examples of unusual behavior of $mm(T^{n})$.

{\bf Theorem.}  {\it One can find an infinite set $A$ of integers and
a weakly mixing automorphism $R$ such that $mm (R^n)=n$ and $mm(R^{n+1}) =1$
for all $n\in A$. Moreover, the cardinalities of the multiplicity sets for~$R^n$
are not bounded. There is a weakly mixing automorphism~$S$ such that
all the powers $S^{n}$ have homogeneous spectrum and the set of limit points
of the sequence $\{\frac{mm (R^n)}{n} \}$ is infinite. }
\vspace{3mm}  

We {\it a posteriori} come to the following principle:  {\it  
given ergodic automorphism $H$ with discrete spectrum
there is  a weakly mixing automorphism  $R$  such that $mm (R^n)= mm (H^n)$ for all $n>0$. }
We apply Ageev's approach based on generic group actions and we
control the spectrum via weak limits of powers. 
The latter 
has an old history  (A.Katok[K], V.I.Oseledets[Os], A.M.Stepin [S] et al).
This note mostly uses [R], [Ry].

We  present some  results  obtained by the author  during his visit  
at {\it Laboratoire de Mathematiques Raphael Salem }  of  Rouen University. 
We thank  el Houcein el Abdalaoui and Thierry de la Rue  for the hospitality
and their interest in this  work. Special thanks go to Sergey Tikhonov for
useful criticism and Jean-Paul Thouvenot for stimulating questions.

\section{  Group actions, homogeneous spectra,  little joinings
}  

Consider a measure-preserving free action of the group $G_p$ generated by
two elements $s$,$\phi$ and the following relations:

1. $\phi^p=e$,

2. the elements $s, \ \phi s\phi^{-1},\ \phi^{2} s\phi^{-2} , \dots , 
\phi^{p-1} s\phi^{-p+1}$ commute.
\\   Denoting the corresponding elements of the action by $S$ and $\Phi$, we now put
$$R:=S\Phi, \quad T:= R^p.$$
It is proved in [Ry] that the automorphism $R$ has
simple spectrum for a generic $G_p$-action. The automorphism~$T$ and all of
its positive powers have homogeneous spectrum of 
multiplicity $p$. The proof of Theorem~1.1 in [Ry] actually shows 
that all the automorphisms $R^n$ have simple spectrum for $n \perp p$ (no common divisors). 
We must only check that Model 2 [Ry]  guarantees the simplicity of spectrum of~$R^n$.
The automorphism $S$ also has simple spectrum for a generic group action.

{\bf  Factors and joinings.}
Our $T$ has a factor $\tilde T$. Indeed,  the automorphism  $T$ commutes with periodic $\Phi$,
we have an $T$-invariant algebra $\cal F$ of $\Phi$-invariant sets.
This factor has simple spectrum [Ry], but  $\F \neq R\F $ (or $\F \neq S\F$).
So there is a non-trivial self-joining of this factor. 
The correspondence between joinings and factors is well known (see [R1], [Th] for
joinings and their relations to factors).  We put
$$ \nu (A\times B) =\int_X  \chi_A \chi_{RA} d\mu $$
for $A,B\in \cal F$.
Since $R$ commutes with $T$, we have
$$\nu (\tilde TA\times \tilde TB) =\int_X  T\chi_A T\chi_{RA} d\mu =\int_X  \chi_{\tilde A} \chi_{R\tilde A} d\mu =
\nu (A\times B).$$
Thus we have self-joinings of $\tilde T$, but the measure $\nu$ is well defined for all measurable sets $A,B$.  
In the theory of joinings, this measure is said to be off-diagonal.  The action of
$T\times T$ with this invariant measure is naturally
isomorphic to $T$ via the map $x \mapsto (x,Rx)$. We take $p=3$ and
consider the map $ x \mapsto (Rx, R\Phi x, R\Phi^2 x)$. There is a set $D$ such that 
$X$ is a disjoint union of $D,\Phi D, \Phi^2 D$.  Our factor $\tilde T$ is 
associated with the map 
$$ \tilde T: D\to D, \quad \tilde Tx =\Phi^{n(x)}Tx,$$
where $n(x)$ is defined by the condition $\Phi^{n(x)}Tx\in D$. We consider the map 
$$ x \mapsto (\Phi^{n_1(x)}Rx, \Phi^{n_2(x)}R\Phi x, \Phi^{n_3(x)}R\Phi^2 x)\in D\times 
D\times D.$$
Since the relation $\Phi^{n}R\Phi^{m}=R$ implies in our group that $n=m=0 (mod p)$, 
we see that the  points $\Phi^{n_1(x)}Rx$, $\Phi^{n_2(x)}R\Phi x$, 
$\Phi^{n_3(x)}R\Phi^2 x$ must be different for almost all $x\in D$.
(Here we use the obvious fact that different automorphisms commuting with 
an ergodic action have pairwise disjoint graphs.) 
So we obtain a bijection between the off-diagonal in $X\times X$ and the
graph of this 3-valued map on~$D$. Thus our joining sits on the graph of 
a 3-valued map. 

{\bf Remarks. 1.} If $\tilde T$ has simple spectrum and  
commutes with an $n$-valued map, then the adjoint map is also $n$-valued.
Indeed, the corresponding commuting Markov operators $J$ and $J^*$ satisfy
$$ \frac{1}{m}I + \frac{m-1}{m}Q = J^*J=JJ^*=\frac{1}{n}I + \frac{n-1}{n}Q', $$
where the graphs of $Q$ and $Q'$ are disjoint from the diagonal. Thus,  
$\frac{1}{m}I= \frac{1}{n}I$,  $m=n$. 

{\bf 2.}  For a generic $G_p$-action, the corresponding $\tilde T$ is of rank one.
Moreover, it is a so-called  flat-stack rank-one automorphism. J.\,King proved 
that each ergodic joining of such an automorphism is a weak limit of
off-diagonal measures ( $\tilde T^{n_i}\to^w J$).

{\bf 3.} S.\,Tikhonov used $G_p$-actions to prove the  existence of a mixing
automorphism~$T$ 
with $p$-homogeneous spectrum.  Thus he also found a non-trivial self-joining
for the mixing
automorphism~$\tilde T$. However, when $p>2$, our factor  $\tilde T$ cannot be 
of rank one because of King's another well-known theorem saying that mixing 
rank-one automorphisms possess the property of minimal self-joinings (see [R1] 
for King's theorem).  So it follows from $\F \neq R\F $ that  
these  factors are independent, hence,  the joining $\nu$ is a product measure.  This contradicts
 the fact that $\nu$ is situated  on a graph. 
In fact the local rank of mixing automorphism $T$ is not greater than $p^{-1}$ 
(see  [Th]).

Problem: {\it must the local ranks of mixing automorphisms $R, \tilde T$ in question be
 zero? } 
We conjecture that our present spectral results can be established for mixing 
automorphisms via Tikhonov's approach [T].

{\bf 4.} It is not hard to obtain a finite-valued self-joining of a mixing map 
with simple spectrum. Suppose that the mixing symmetric product~$T^{\odot 3}$ 
has simple spectrum (see~[R]). Consider $T\otimes T\otimes\ T$ and two factors:
the symmetric factor $\cal F$ and a perturbation of it, say, $(I\otimes I\otimes\ T)\F$. 
These two factors generate  the  algebra of all measurable sets. 
The corresponding  joining of $T^{\odot 3}$ will be finite-valued. 
The dynamics on this joining exactly imitates~$T\otimes T\otimes\ T$.
We also note that  $mm(T\otimes T\otimes\ T)=3!$ and the spectrum of 
$T\otimes T\otimes\ T$ is non-homogeneous (type~(3,6)).

{  \bf 5. Hard problems}. 
Find \\
 5.1.  an ergodic  $T$ such that     $mm(T^{p^n})=1$,  $mm(T^{p^{n+1}})=p$  for some  $p,n>1$;
\\
5.2.   an ergodic automorphism $T$ such that     $mm(T^{n})$ is not bounded, but  $mm(T^n)< ln ln (n)  + Const$;
\\
5.3.  an  automorphism $T$ such that     $Rank (T^{p^n})=1$,  $Rank(T^{p^{n+1}})=p$ for some  $p,n>1$;
\\
5.4.  a flow  $T_t$ such that     $mm(T_s)=1$, $0<s<1$,  $mm(T_2)=2$  (and variations of this problem including rank analogies).

{ \bf 6. \it Easy problem}. Find a weakly mixing flow  $T_t$ such that
$mm(T_r)=\infty $ for all rational $r$ and $mm(T_a)=1$ for all irrational~$a$. 
We note that there is a discrete-spectrum model, so the reader could apply
our principles and methods.

\section{  Maximal spectral multiplicity for the powers. Periodic cases}  

In [R] we edowed a weakly mixing automorphism $T$ with ($G_p$-generic) properties:
arbitrary polynomial weak limits for its powers. The   results of [R],[Ry]
provided  all the powers  $T^N$, $N > 0$,  homogeneous spectrum of multiplicity $p$.
We note that [A1] contains non-weakly-mixing  examples 
with $mm(T^n)$  of the type~(1,2,1,2,...).  Weakly mixing and mixing examples were
found  by the author (see~[R]).  Here we
consider $R(x,y)=(y,Tx)$ 
on~$X\times X$.  It is possible to have  $mm(R^{2n})= mm(T^n \otimes T^n)=2$ and
$mm(R^{2n+1})=1$.  
\vspace{4mm} 

{ \bf  The  case  $\bf \{1,1,\dots,1, p,1,1,\dots,1,p,1,\dots\}$.}
We look at $mm(R^{n})$ for $R$ in a generic $G_p$-action.
Let $p$ be a prime. Then we obtain $mm(R^n)=1$
for all $n\perp p$ (no common divisors). We have $mm(R^{np})=p$.
\vspace{4mm} 

{ \bf The case  $\bf \{1,\dots,1, p,1,\dots,1,q,1,\dots,1,pq, 1,\dots\}$.}
Let $p,q$ be primes. Then the desired behavior of $mm(R^{pqN})$ is easily
seen to hold for the $G_{pq}$-action with $mm(R^{pqN})=pq$.
However we propose a more complicated solution to explain our strategy 
for the infinite (inverse) limit construction.

We fix the primes $p,q$ and consider the transformations $R_1$, $R_2$ for 
the actions of $G_p$, $G_q$ respectively. Put $T_1=R_1^{p}$, $T_2=R_2^{q}$. 
Our first aim is to find conditions that guarantee the simplicity of spectrum
of $R=R_1\times R_2$ and, moreover, we wish to have

$mm(R^n)=pq$  if   $n=n'pq$, $n'>0$;

$mm(R^n)=1$ if   $n\perp pq$;

$mm(R^n)=p$  if   $n=n'p$ and $n'\perp q$;

$mm(R^n)=q$  if   $n=n'q$ and $n'\perp p$;
\\
To obtain this spectral behavior of the powers, we use the following weak limit
for some sequence $n_i$ (condition WL(2)):

$$(R_1\otimes R_2)^{pn_i} \to I\otimes R_2.$$
\\
This implies for any fixed integer $N>0$ that
$$(R_1\otimes R_2)^{Npn_i} \to I\otimes R_2^N.$$
\vspace{4mm} 

{\bf LEMMA 1.}(See Lemma 3.2 in [Ry]) {\it  Suppose that $R$ and $R'$ have simple
spectra. If there is a sequence $n_i$ such that
$(R\otimes R')^{n_i} \to I\otimes R'$, then $R\otimes R'$ has  simple spectrum.}
\vspace{4mm}

So for all $N\perp pq$ we obtain from WL(2) and Lemma~1 that  
$$mm(R_1^N\otimes R_2^N)=1$$
since $R_1^N$ and $R_2^N$ have simple spectra. 
So, for example, we get
$$mm(R)=mm(R^{pq-1})= mm(R^{pq+1})=1$$
and
$$mm(R^{pq})=pq.$$
The last is a direct consequence of the following facts:

1. $R^{pq}=T_1^q\otimes T_2^p$ is a product of two automorphisms with multiplicities
$p$ and $q$, whence $mm(R^{pq})\geq pq$;

2. $mm(R^{pq})\leq pq$ since $mm(R)=1$.
\\
Combining these two facts, we get $mm(R^{pq})=pq$.

Thus we must explain how to achieve~WL(2).  
We find a sequence $n_i\to\infty$ and a related sequence $\tilde{n}_i$ such that
$$pn_i = q\tilde{n}_i +1.$$
As shown in [Ry], a generic $G_p$-action has the property
$$T_1^{n_i}\to I$$
for some subsequence of indices $i$.
But for a generic $G_q$-action we have
$$T_2^{\tilde n_i}\to I$$
for some subsubsequence of indices $i$.
  Denoting this subsubsequence again by $n_i$, we get WL(2). Indeed,
$$R_1^{pn_i} = T_1^{n_i} \to I,$$
$$  R_2^{pn_i} =R_2^{q\tilde{n}_i +1}=T_2^{\tilde{n}_i}R_2\to R_2.$$ 

\section{  Conditions WL(k) imply that $mm(R^{n})=n$    and 
$mm(R^{n+1})=1$  for all $n$ of the form $n=p_1p_2\dots p_k$}

We define condition WL(3) for a triple of   $G_{p_1}, G_{p_2}, G_{p_3}$-actions
as the existence of  sequences $n_i(1)$ and  $n_i(2)$ such that

$$(R_1\otimes R_2)^{p_1n_i(1)}\to I \otimes R_2,$$

$$(R_1\otimes R_2\otimes R_3)^{p_1p_2n_i(2)}\to I \otimes I\otimes R_3.$$

Definition.  Property WL(k) for a $k$-tuple of $G_{p_1}, \dots, G_{p_k}$-actions
means that there are sequences $n_i(1)$,   $\dots, n_i(k)$ such that

$$(R_1\otimes R_2)^{p_1n_i(1)}\to I \otimes R_2,$$
$$\dots$$
$$(R_1\otimes \dots\otimes R_k)^{p_1\dots p_{k-1}n_i(k-1)}\to I \otimes\dots I\otimes R_k.$$

If all operators $R_1, R_2,\dots$ have simple spectra and satisfy conditions 
WL(k), $k>1$, then we see by induction that their finite direct products have
simple spectra. This yields
$$mm(R_1\otimes R_2\otimes\dots)=1.$$

Put $R=R_1\otimes R_2\otimes\dots$ and suppose that $N\perp p_1,p_2,\dots$. 
We claim that $$mm(R^N)=1.$$  Indeed, replacing the sequence $n_i(k)$ 
in conditions WL(k+1) by $Nn_i(k)$, we obtain conditions WL(k+1) for the 
automorphisms  $R_1^N, R_2^N,\dots$ which also have simple spectra.

We now  explain why  $$mm(R^{p_1\dots p_k}) = p_1\dots p_k.$$
Since $mm(R)=1$, we know that
$$mm(R^{p_1\dots p_k}) \leq p_1\dots p_k.$$   

Let $$q(m)=\frac {p_1p_2p_3\dots p_k}{p_m}.$$
We know that  $$mm(T_m^{q(m)})=p_m.$$
Thus the representation
$$R
^{p_1\dots p_k}= T_1^{q(1)} \otimes T_2^{q(2)}\otimes\dots\otimes T_k^{q(k)}
\otimes Q_k$$ yields that
$$mm(R^{p_1\dots p_k})\geq p_1\dots p_k,$$   
whence 
$$mm(R^{p_1\dots p_k})= p_1\dots p_k.$$  

\section { How to obtain conditions WL(k)}

{\bf LEMMA 2.} {\it  Let $\tilde A$ and $A$ be infinite sets of integers.
Then for every $p>0$ there are two sequences
$\tilde n_i(2)\in \tilde A$, $\tilde n_i(2)\in A$ and an automorphism $T$
from a generic action of~$G_p$ 
such that   $T^{\tilde n_i}\to I$  and $T^{ n_i}\to I$.}
\vspace{4mm}
\\
{\bf REMARK 3.} {\it We can actually find the convergents $T^{\tilde n_i}\to P(T)$
and $T^{ n_i}\to Q(T)$ for arbitrary admissible polynomials $P,Q$.
(A polynomial $P(T)=\sum a_iT^i$ is said to be admissible if $a_i\geq 0$ 
and $\sum a_i=1$.)
}
\vspace{4mm} 

Using Lemma 2 we find a rigid sequence $n_i(1)$
for $T_1$:  $T_1^{n_i(1)}\to I$.  Then we find $T_2$ with rigid sequences $\tilde n_i(2)$ and 
$ n_i(2)\subset n_i(1)$, where
$ n_i(2)$ is a subsequence of  $n_i(1)$ and the sequence $\tilde n_i(2)$
is related to $ n_i(1)$ by the equation 
$$p_1n_i(2) = p_2\tilde{n}_i(2) +1.$$
We now have
$$(R_1\otimes R_2)^{p_1n_i(2)}=(T_1^{n_i(2)}\otimes T_2^{\tilde n_i(2)}R_2)\to  I \otimes R_2.$$

Then we find $T_3$ with rigid sequences $\tilde n_i(3)$ and $ n_i(3)\subset n_i(2)$  such that
 $\tilde n_i(3)$  is related to $ n_i(3)$ by the equation 
$p_1p_2n_i(3) = p_3\tilde{n}_i(3) +1$.
Thus we obtain
$$(R_1\otimes R_2\otimes R_3)^{p_1p_2n_i(3)}=(T_1^{p_2n_i(3)}\otimes T_2^{p_1n_i(3)}\otimes T_3^{\tilde n_i(3)}R_3)\to 
 (I \otimes I \otimes R_3).$$

Using induction, we find by Lemma 2 an automorphism $T_k$ with rigid sequences 
$\tilde n_i(k)$ and $ n_i(k)\subset n_i(k-1)$  such that
 $\tilde n_i(k)$  is related to $ n_i(k)$ by the equation 
$$p_1p_2\dots p_{k-1}n_i(k) = p_k\tilde{n}_i(k) +1.$$
We get $T_1,T_2, \dots $ for which the corresponding set $R_1,R_2, \dots $ has 
properties WL(k).

Thus the product $R=R_1\otimes R_2\otimes R_3\otimes\dots $ has   simple spectrum.
Using section 3, we obtain 

\vspace{4mm}
{\bf THEOREM 4.} {\it  Let $P=\{p_1, p_2,p_3,\dots \}$, where   $p_1< p_2<p_3<\dots $ are primes. 
Suppose that $N=p_{k_1}^{d_1}  p_{k_2}^{d_2}\dots p_{k_m}^{d_m}q$,
$d_j>0$, and $q\perp P$. Then $mm(R^N)= p_{k_1}  p_{k_2}\dots p_{k_m}$ and $mm(R^q)=1$.
For example, if   $p_{k+1} > p_1p_2\dots p_k$, then    $mm(R^{p_1p_2\dots p_k +1})=1$.}
\vspace{4mm}

We note that the above results were partially given in [Ry] 
for finite products. The WL(k)-methods were used there to show that
$mm(R^N)= p_{k_1}  p_{k_2}\dots p_{k_m}$ for $N= p_{k_1}  p_{k_2}\dots p_{k_m}M$
(the similar methods were also developed in the proof of Theorem 3.3 in [Ry]).
We just applied this approach to infinite products.  

Our next aim is to pass to homogeneous spectra for the powers. Before doing this,
we note that the cardinality $cardm(R^N)$ of the set of spectral multiplicities for 
$R^N$ can also be controlled by the same methods.
\vspace{4mm}

{\bf THEOREM 5.} {\it Suppose that $N=p_{k_1}^{d_1}p_{k_2}^{d_2}\dots p_{k_m}^{d_m}q$,
$d_j>0$, and $q\perp P$. Then $cardm(R^N)= 2^m $ and $cardm(R^q)=1$. For example, the
 multiplicity set of $R^{p_ip_j}$  is  $\{1,p_i,p_j, p_ip_j\}$.
}
\vspace{4mm}

The proof is a natural adaptation of the arguments proving
Theorem 3.3 in [Ry] to the case of infinite products.

\section{Homogeneous  spectral multiplicities  of the powers}

Let $R_p$ denote  a cyclic permutation on the space $\{1,2,\dots,p\}$ with 
the  uniform measure. We consider
$$R=R_{p_1}\times R_{p_2}\times R_{p_3}\times\dots,$$
where ${p_i}\in P$ are different prime numbers. Our measure space $(X,\mu )$ 
is the group
$$X=\Z_{p_1}\times \Z_{p_2}\times \Z_{p_3}\times\dots$$
with Haar measure $\mu$.
The transformation $R$ has rank one, so it is ergodic with simple discrete spectrum.
The same holds for the powers $R^n$ with any $n\perp P$.
But if $n=p_1p_2\dots p_k$, then we have $mm(R^n)=n$.

We consider  the following group of automorphisms of the group $X$: 
$$A=\Z_{p_1-1}\otimes \Z_{p_2-1}\otimes \Z_{p_3-1}\otimes\dots.$$
Suppose
that $q\in Q\perp P$. Then $R^q$ is conjugate to~$R$. 
For every $q\in Q \perp P$ we find $\Psi_q\in A$ such that
$$\Psi_q^{-1}R\Psi_q = R^q.$$

We define a group A-R generated by $R$ and all $\Psi_q$, $q\in Q$.  
Using some ideas of [A],[D], we easily obtain the following result.
\vspace{4mm}

{\bf Theorem 6.} {\it Suppose that $Q \perp P$ are disjoint subsets  of prime numbers
and $P$ is infinite. Then there is a rank-one weakly mixing transformation~$R$ 
conjugate to all the powers $R^n$ with $n\in Q$}. 
\vspace{4mm}

{\it Proof}.  The  above A-R-model (with a rank-one automorphism $R$ of discrete spectrum) 
is amenable and free. We add a consideration of the Bernoulli action to
obtain a weakly mixing A-R-generic automorphism~$R$. The  rank-one property is generic. 
We automatically obtain the conjugations in question. The proof is complete.  

We return to $R=R_{p_1}\times R_{p_2}\times R_{p_3}\times\dots$. 
The power $R^n$ is non-ergodic for $n=p_iN$. 
But all powers have homogeneous spectra. This suggests  
the following theorem, where we write $hm(T)=n$ if and only if $T$ has
homogeneous spectrum of multiplicity~$n$.

\vspace{3mm}
{\bf THEOREM 7.} {\it Let  $P=\{p_1, p_2, \dots\}$ be an infinite set of 
prime numbers. Then there is a weakly mixing automorphism $R$ such that   
we have $hm(R^N)=p_{k_1}p_{k_2}\dots p_{k_m}$ and $hm(R^q)=1$
for all  $N=p_{k_1}^{d_1}  p_{k_2}^{d_2}\dots p_{k_m}^{d_m}q$,  $d_j>0$,  
and all $q\perp P$.}

\vspace{3mm} 
{\it Proof}. We slightly generalize some ideas of [Ry]. Instead of the 
$\Z_p$-subaction, we consider a subaction of the group     
$$\Z_{p_{1}}\oplus \Z_{p_{2}}\oplus \Z_{p_{3}}\oplus\dots.$$
Let $G_P$ be the group generated by the elements
$R,F_1, F_2, F_3 \dots$ and the following relations:

1. $F_iF_j=F_jF_i$, \  $F_i^{p_i}=I$,

2.  the elements $F_i^{-n}RF_i^{n-1}$ commute.

\vspace{3mm}
{\bf LEMMA 8.} {\it   A generic $G_P$-action possesses the following properties:

1. $R^q$ has a simple continuous spectrum if $q\perp P$.

2. $F_i^{-1}R$ has a simple spectrum.

3. The  action of any infinite order group element  have continuous spectrum.

4.  $R^{p_i}$ has  homogeneous spectrum of multiplicity $p_i$. 

5. $Rank (R^n)=1$ for $n\perp P$, and $Rank(R^n)=mm(R^n)$.}
\vspace{3mm}

{\it Proof of Lemma 8}.  

{\bf Model 1.} Property 1 can be realized using our representation 
$R=R_1\otimes R_2\otimes R_3\otimes\dots $ from section 5.
We only note that the product
$R_1\otimes R_2\otimes R_3\otimes\dots $ is simply connected
and admits a natural free $G_P$-action. Indeed, put
$$F_i = I\otimes I \dots\otimes \Phi_i  \otimes I\dots,$$
where $\Phi_i$ is taken from the $G_{p_i}$-action.  Since the
$G_{p_i}$-action is free, so is the $G_{P_i}$-action.
\vspace{2mm} 

{\bf Model 2.} To realize property 2, we consider $R_i$ acting 
on $(X_i,\mu_i)=(\tilde X^{p_i},\tilde \mu^{p_i})$ by the formula
$$R_i =F_i\circ (U_1\otimes U_2\otimes\dots   U_{p_i}).$$
We achieve a simple spectrum for such $R_i$, say, taking $U_k$ in a generic rank-one
flow. Then
$$F_i^{-1}R= R_1\otimes\dots R_{p_{i-1}}\otimes (U_1\otimes U_2\otimes\dots
U_{p_i})\otimes R_{p_{i+1}}\otimes\dots .$$
Using the WL(k)-method, we again achieve a simple spectrum for this product.  
Thus the automorphism $R$ and all cyclic coordinate permutations $F_i$ 
generate a free $G_P$-action with property~2.

The Bernoulli action on Model 2 gives us property 3. 

A generic $G_P$-action has all generic properties 1,2,3.  They together imply
property 4.  Indeed, consider
the $\Z^{p_i}$-subaction generated by the commuting automorphisms  

 $S_1=F_i^{-1}R$,\  $S_2=F_i^{-2}RF_i^{1}$,    $\dots$,  $S_{p_i-1}=F_i^{{-p_i}+1}RF_i^{{p_i}-2}$, \ $S_{p_i}=RF_i.$

The spectrum of this action is simple since the spectrum of $S_1$ is simple. 
We repeat the proof of Lemma 1.2 [Ry] almost word-by-word and
obtain that the automorphism $S_1S_2\dots S_{p_i}=R^{p_i}$ has
homogeneous spectrum of multiplicity ${p_i}$.  The lemma is proved.

Theorem 7 follows from Lemma 8, which enables us
to replace $mm(R^n)$ by $hm(R^n)$ in Theorem~4.  
The author plans to represent more details in a journal version.

{\bf Concluding remarks. 1.}
 It remains to find $R$ with   the following properties:
$hm(R^{n}) =n, {\  \ }  n=1,2,3, \dots, 2009, 2010,   {\ } {\ }
hm(T^{2011}) =1.$
To do this, we consider instead of the primes $p_{k_1},  p_{k_2},\dots$
their powers
$p_{k_1}^{d_1},  p_{k_2}^{d_2},\dots$
and use the following fact. If $hm(T)=1$ and $hm(T^{p^d})=p^d$, 
then $hm(T^{p^m})=p^m$ for all $m<d$.

Of course, the author was lucky since 2011 is prime.
So we consider all prime numbers and put $d_i=1$ for $p_i=2011$,  
and $d_j=10$ for all $p_j\neq 2011$.
Choosing the space
$$X=\Z_{2^{10}}\times \Z_{3^{10}}\times \Z_{5^{10}}\times\dots
\times \Z_{2003^{10}}\times \Z_{2011}\times \Z_{2017^{10}}\times\dots,
$$
we obtain required $R$ with discrete spectrum and, via an effort of will
as above, we transform the discrete spectrum into  a continuous one.
\vspace{4mm} 

{\bf 2.} There is a simpler way to get
$$hm(R^{n}) =n, {\  \ }  n=1,2,3, \dots, 2009, 2010, {\ } {\ }
hm(T^{2011}) =1.$$
Consider a generic $G_n$-action with $n=2010!$ (see [Ry]).
The corresponding $R$ has the property above, but the behavior of 
$hm(R^{n})$ is periodic. However, Tikhonov's approach enables us
to make $R$ mixing.

{\bf 3. Connected problems.}   
(a) Let  ${\cal A}_{G_{n!}}$ denote   the space of  $G_{n!}$-actions.   We  consider the  projections 
$\pi_n:{\cal A}_{G_{(n+1)!}}\to{\cal A}_{G_{n!}}$  defined by
$\pi_n \Phi =\Phi^n,$ $\pi_n R =R^n. $
Prove  that  for all $n$  the image of $\pi_n$  is residual in ${\cal A}_{G_{n!}}$.

(b)  Find   mixing automorphisms with the following  properties:

1.  $\{hm(R^{n}) : n\in\N\}  =\{1,2,3,4,  \dots,k, k+1, \dots \}$;

2.   $\{hm(R^{n}) : n\in\N\}  =\{1,  p,  p^2, p^3,  \dots, p^k, \dots\}$.

\large
\begin{center}
{\bf References}
\end{center}

{\bf [Ab]}  E. H. El Abdalaoui. On the spectrum of the powers of Ornstein transformations,  Sankhya, Ser. A 62 (2000), no 3, 291 -306.

{\bf  [Ag] }   O. N. Ageev, Dynamical systems with an even-multiplicity Lebesgue component in the spectrum.
 Mat. Sb. 178:3 (1988), 307-319

 {\bf [Ag1]} O. N. Ageev, The spectral type of the rearrangements $T_{\alpha, \beta}$. Mat. Sb., 188:8 (1997), 13-44

{\bf [Ag2]} O. N. Ageev. Spectral rigidity of group actions: applications to the case gr$(t, s; ts = st^2)$.
Proc. Amer. Math. Soc. 134 (2006), 1331-1338.

{\bf [D]}  A. I. Danilenko. Weakly mixing rank-one transformations conjugate to their squares.    
    Studia Math. 187 (2008), 75-93

{\bf [K]}  
A. Katok, Combinatorial constructions in ergodic theory and dynamics, Univ. Lecture Ser.,
vol. 30, Amer. Math. Soc., Providence, RI 2003.

{\bf [L]} M. Lemanczyk, Toeplitz $Z_2$-extensions, Ann. Inst. H. Poincar\' e 24 (1988), 1-43.

{\bf [Os]} V. I. Oseledets, The spectrum of ergodic automorphisms. Dokl. Akad. Nauk SSSR, 168:5 (1966), 1009-1011 

{\bf [R1]}  V.V.  Ryzhikov.  Mixing,   rank, and minimal self-joining of actions with an invariant measure. 
 Mat. Sb., 183:3 (1992), 133-160 

{\bf [R]} V.V. Ryzhikov. Weak limits of powers, simple spectrum of symmetric products, and
rank-one mixing constructions.  Sb. Math. 198:5 (2007), 733-754.

{\bf [Ry]} V.V. Ryzhikov. Spectral multiplicities and asymptotic operator
properties of actions with invariant measure. Sb. Math., 200:12, (2009), 1833-1845 

{\bf [S]} A. M. Stepin, Spectral properties of generic dynamical systems, Izv. Akad. Nauk SSSR
Ser. Mat. 50:4 (1986), 801-834.

{\bf [Th]}   J.-P. Thouvenot. Some properties and applications of joinings in ergodic theory. Ergodic theory and its connections with harmonic analysis, Proc. of the 1993 Alexandria Conference, LMS Lecture notes series, 205, Cambridge Univ. Press, Cambridge, 1995.  

{\bf [T]} S.V. Tikhonov.  On mixing transformations with homogeneous spectrum.  Sb. Math. (to appear).
\vspace{10mm}

Department of Mechanics and Mathematics, 

Lomonosov Moscow State University,

GSP-1, Leninskie Gory, Moscow, 119991,

 Russian Federation

E-mail address: vryzh@mail.ru

\end{document}